\documentclass[preprint,12pt]{elsarticle}

\usepackage{amssymb}
\usepackage{amsmath}
\usepackage{amsthm}

\newtheorem{remark}{Remark}
\newtheorem{theorem}{Theorem}[section]
\newtheorem{defination}{Definition}[section]

\newtheorem{lemma}[theorem]{Lemma}
\usepackage{color,amsmath}

\usepackage{subfig}

\journal{Applied Numerical Mathematics}

\begin{document}

\begin{frontmatter}



\title{Jacobi convolution series for Petrov-Galerkin scheme and general fractional calculus of arbitrary order over finite interval} 


\author[aff]{Pavan Pranjivan Mehta}
\ead{pavan.mehta@sissa.it}
\author[aff]{Gianluigi Rozza\corref{CorrespondingAuthor}} 
\ead{grozza@sissa.it}
\cortext[CorrespondingAuthor]{Corresponding author}
\affiliation[aff]{organization={Mathlab, Mathematica area, 
SISSA, International School of Advanced Studies},
            addressline={Via Bonomea, 265}, 
            city={Trieste},
            postcode={34136}, 
            state={TS},
            country={Italy}}

\begin{abstract}
Recently, general fractional calculus was introduced by \citet{kochubei2011general} and \citet{luchko_sonin_21} as a further generalisation of fractional calculus, where the derivative and integral operator admits arbitrary kernel. Such a formalism will have many applications in physics and engineering, since the kernel is no longer restricted. We first extend the work of \citet{luchko_finite_interval} on finite interval to arbitrary orders. Followed by, developing an efficient Petrov-Galerkin scheme by introducing Jacobi convolution series as basis functions. A notable property of this basis function, the general fractional derivative of Jacobi convolution series is a shifted Jacobi polynomial. Thus, with a suitable test function it results in diagonal stiffness matrix, hence, the efficiency in implementation. Furthermore, our method is constructed for any arbitrary kernel including that of fractional operator, since, its a special case of general fractional operator. 
\end{abstract}



\begin{keyword}
General fractional calculus \sep Jacobi convolution series  \sep Petrov-Galerkin scheme \sep spectral methods




\MSC[2020] 26A33 \sep 35R11 \sep 65M70 \sep 65N35 \sep 33C45 \sep 42C05

\end{keyword}

\end{frontmatter}



\section{Introduction} \label{intro}

Fractional calculus is a natural extension of standard integer-order calculus. The primary aim of such a generalisation is to extend the notion of derivatives and integral with orders defined in $\mathbb{R_+}$. Paradoxically, for non-integer orders, these operators are nonlocal. Inorder to extend integer-order operator to fractional operators, two approaches exist: 

\begin{itemize}
    \item Starting from the limit definition of derivatives, we derive to Gr\"{u}nwald-Letnikov derivative (\ref{eq:1.1}) \cite{grunwald1867uber} (see also \cite{podlubny1999})

\begin{equation} \label{eq:1.1}
        f^{(p)}(x) = \frac{d^{p} f} {d x^{p}} = \lim_{h \rightarrow 0} \frac{1} {h^{p}} \sum_{r=0}^{N} (-1)^{r} \binom{p}{r} f (x - r h)
\end{equation}
    
    \item The second direction involves generalising the Cauchy repeated integration formula. Recall, Cauchy repeated integration formula (\ref{eq:2.6}) for $p \in  \mathbb{N}$, 

\begin{equation}  \label{eq:2.6}
   \underbrace{ \int_{a}^{x} \int_{a}^{x_p} \int_{a}^{x_{p-1}} ... \int_{a}^{x_{2}}}_{\text{p-integrals}} f(x_1) d x_1  ... d x_p = {1 \over (p-1)!}  \int_{a}^{x} {(x - \tau)}^{p-1} f (\tau) d \tau
\end{equation}

    We now invoke, $\Gamma(p) = (p-1)!$, where $\Gamma(.)$ is the Euler gamma function, thus the Riemann-Liouville fractional integral is defined as (\ref{eq:RLI}) for ($p \in \mathbb{R}_+$)

\begin{equation}  \label{eq:RLI}
   {}_a I^p_x f(x)  ~:=~ {1 \over \Gamma(p)}   \int_{a}^{x} {(x - \tau)}^{p-1} f (\tau) d \tau
\end{equation}
    
    This leads to the definition of Riemann-Liouville derivative as (\ref{eq:RLD}) \cite{podlubny1999}.

\begin{equation}  \label{eq:RLD}
{}_a^{RL} D^{(p)}_x f(x)  ~:=~ {1 \over \Gamma(n-p)} {d^n \over d x^n}
\int_{a}^{x} {(x - \tau)}^{n-p-1} f (\tau) d \tau ,  n-1 \leq p < n 
\end{equation}


\item  Besides, these two definitions, Caputo's definition \cite{Caputo1967} for derivative is given as (\ref{eq:CD}) \cite{podlubny1999}. 

\begin{equation}  \label{eq:CD}
{}_a^{C} D^p_x f(x)  ~:=~ {1 \over \Gamma(n-p)} 
\int_{a}^{x} {(x - \tau)}^{n-p-1}   f^{(n)} (\tau) d \tau ~~, ~ n-1 < p \leq  n 
\end{equation}

\end{itemize}

\noindent
It is to be noted that, Caputo derivative mitigates two key problems (sec. 2.4 of \cite{podlubny1999}) faced by Riemann-Liouville derivative, 

\begin{itemize}
    \item Caputo derivative of  constant function is zero, whereas for Riemann-Liouville derivative is generally not true. 
    
    \item Initial conditions for Caputo derivative are prescribed in a classical sense as opposed to Riemann-Liouville derivative (sec. 2.4 of \cite{podlubny1999}).
\end{itemize}

Owing to these two facts, the Caputo derivative is often used for applications such as turbulence \cite{mehta2019discovering, mehta2023fractional}. It is important to note that, both the Riemann-Liouville and Caputo derivative satisfy both first and second fundamental theorem of calculus  \cite{podlubny1999, diethelmbook}, akin to classical integer-order calculus.

An overwhelming question arises for physicists and engineers, \textit{can a operator with a power-law kernel describe all physical processes?} The answer is clearly no, and there exists examples in turbulence studies, where the kernel is found to be other than a power-law \cite{hamba2017history} (following \cite{kraichnan1987eddy}). However, the mathematical theory of such generalisation is recent \cite{kochubei2011general, luchko2021general}.


With regards to generalisation of fractional operators, Sonine \cite{sonine1884generalisation} recognised a key property that the convolution of the kernel (of the fractional derivative and integral) is unity, thereafter proposed the condition (\ref{eq:sonin}) for any pair ($k(x), \kappa(x)$) for analytical solution. 

\begin{equation} \label{eq:sonin}
    \int_0^x k (x-t) \kappa (t) dt = 1, ~x > 0.
\end{equation}

Indeed, there are more than one examples of kernels satisfying  (\ref{eq:sonin}). For instance, Sonine introduced the kernel pair in \cite{sonine1884generalisation} of the form (\ref{eq:sonine_par}) (following the notations of \cite{luchko2021general}). 

\begin{align} \label{eq:sonine_par}
    \begin{split}
        &\kappa(t) = h_\alpha (t) . \kappa_1 (t), ~~~ \kappa_1 (t) = \sum_{k=0}^{+\infty} a_k t^k, ~ a_0 \neq 0, ~~ 0 , ~~ \alpha \in (0,1) \\
         &k(t) = h_{1-\alpha} (t) . k_1 (t), ~~~ k_1 (t) = \sum_{k=0}^{+\infty} b_k t^k, ~ b_0 \neq 0 
    \end{split}
\end{align}

where, $h_{\alpha}(t) = {t^{\alpha-1}}/{\Gamma(\alpha)}$ and the coefficients follow the relationship (\ref{eq:coef}). For more non-trival examples refer \cite{luchko_finite_interval, luchko2021general, sonine1884generalisation}

\begin{align} \label{eq:coef}
    \begin{split}
        &a_0 b_0 ~=~ 1,~  n = 0 \\
        &\sum_{k=0}^n \Gamma(k+1-\alpha) \Gamma(\alpha + n - k) a_{n-k} b_k = 0,~  n \geq 1.
    \end{split}
\end{align}

\noindent
Perhaps, the first results within the framework of fractional calculus was done in \cite{kochubei2011general} (now known as general fractional calculus), where he introduced a class of kernels which satisfy the following conditions (Kochubei class), 

\begin{itemize}
    \item The Laplace transform of $k$ is $\Tilde{k}$, 
    \begin{equation*}
        \Tilde{k}(p) = (\mathcal{L} k)(p) = \int_0^\infty k(t) e^{-pt} dt
    \end{equation*}
    exists for all $p>0$.

    \item $\Tilde{k}(p)$ is a Stieljes function

    \item $\Tilde{k}(p) \to 0$ and $p \Tilde{k}(p) \to +\infty$ as $p \to +\infty$

    \item $\Tilde{k}(p) \to +\infty$ and $p \Tilde{k}(p) \to 0$ as $p \to 0$
\end{itemize}

However, working with Laplace transform of the kernel is rather cumbersome, thus in \citep{luchko_sonin_21}, another class of kernels were introduced (Luchko class). Followed by its extension to arbitrary order by introducing a modified Sonine condition in \cite{luchko2021general} and finite interval in \cite{luchko_finite_interval}. The work on finite interval for arbitrary orders is an open question, thus in the subsequent section we first extend the general fractional calculus over finite intervals to arbitrary orders. Followed by the development of Petrov-Galerkin scheme for such generalised operator definitions by introducing a new class of basis functions, namely, Jacobi convolution series.

The structure of the paper is as follows:

\begin{itemize}
   \item Section \ref{sec:gfc}: We develop general fractional calculus over finite interval for arbitrary order by extending the work of \cite{luchko_finite_interval}.
    \item Section \ref{sec:basis}: We construct new type of basis functions, namely, the Jacobi convolution series.
    
    \item Section \ref{sec:pg}: We develop a Petrov-Galerkin scheme for general fractional operators.
\end{itemize}

\section{General Fractional Calculus} \label{sec:gfc}

Engineering problems often encounters problems defined on finite intervals, invoking the need for a mathematical theory of general fractional calculus on finite interval. The case for $n=1$ was done in \cite{luchko_finite_interval}. In this section, we shall generalise the results of \cite{luchko_finite_interval} to arbitrary order on finite interval by taking our inspiration from \cite{luchko2021general} for semi-infinite domains. We start by defining the Sonine condition as,

\begin{defination} (see \cite{luchko_finite_interval})
    The pair $(k, \kappa)$ satisfy the left  Sonine condition on an interval $(a, b]$, where $a, b \in \mathbb{R}$, is given by, 

    \begin{equation} \label{eq:left_soncon_1}
        \int_{a}^{x} k (x - t) \kappa (t) dt ~=:~ \{1\}_l,~ a < x \le b.
    \end{equation}
    where $\{1\}_l$ is a function uniformly equal to one over the interval.
\end{defination}

\begin{defination} (see \cite{luchko_finite_interval})
    The pair $(k, \kappa)$ satisfy the right Sonine condition on an interval $[a, b)$, where $a, b \in \mathbb{R}$, is given by, 

    \begin{equation}  \label{eq:right_soncon_1}
        \int_{x}^{b} k (x - t) \kappa (t) dt ~=:~ \{1\}_r,~ a \le x < b.
    \end{equation}
    where $\{1\}_r$ is a function uniformly equal to one over the interval.
\end{defination}

In order to generalise the results of \cite{luchko_finite_interval} for general fractional calculus to arbitrary orders ($n \in \mathbb{N}$) defined over a finite interval, we introduce the modified Sonine condition (introduced in \cite{luchko2021general} for semi-infinite domains), where the kernels ($k_n, \kappa_n$) satisfies the condition as follows,

\begin{defination}
    The pair $(k_n, \kappa_n)$ satisfy the left modified Sonine condition for $n \in \mathbb{N}$ on an interval $(a, b]$, where $a, b \in \mathbb{R}$, is given by, 

    \begin{equation}
        \int_{a}^{x} k_n (x - t) \kappa_n (t) dt ~=~ \frac{ (x - a)^{n-1}} { (n-1)!}~=: \{1\}^{n}_l,~ a < x \le b,~ n \in \mathbb{N}.
    \end{equation}
    where $\{1\}_l$ is a function uniformly equal to one over the interval and
\begin{equation*}
    \{1\}^{n}_l ~:=~  \underbrace{\{1\}_l * \{1\}_l * \dots * \{1\}_l}_{n-terms}.
\end{equation*}
\end{defination}

\begin{defination}
    The pair $(k_n, \kappa_n)$ satisfy the right modified Sonine condition for $n \in \mathbb{N}$ on an interval $[a, b)$, where $a, b \in \mathbb{R}$, is given by, 

    \begin{equation}
        \int_{x}^{b} k_n (x - t) \kappa_n (t) dt ~=~ \frac{ (b - x)^{n-1}} { (n-1)!}~=: \{1\}_r^{n},~ a \le x < b,~ n \in \mathbb{N}.
    \end{equation}
    where $\{1\}_r$ is a function uniformly equal to one over the interval and 
    \begin{equation*}
        \{1\}_r^{n} ~:=~  \underbrace{\{1\}_r * \{1\}_r * \dots * \{1\}_r}_{n-terms}
    \end{equation*}   
\end{defination}

The above formula is a direct consequence of Cauchy repeated integration formula. We follow the convention where, $k_0, \kappa_0$ leads to a zeroth order operator. 

An important example of the kernels, satisfying the left-modified Sonine Condition  is,

\begin{align}
    \begin{split}
         &k_n ~=~ \int_a^x k_{n-1}(x-t) k_1 (t) dt ~=~    k_{n-1} * k_1 ~=~ \underbrace{k_1 * k_1 * \dots * k_1}_{\text{n-terms}}  \\
          &\kappa_n ~=~ \int_a^x \kappa_{n-1}(x-t) \kappa_1 (t) dt ~=~    \kappa_{n-1} * \kappa_1 ~=~ \underbrace{\kappa_1 * \kappa_1 * \dots * \kappa_1}_{\text{n-terms}}  
    \end{split}
\end{align}

where, $k_1, \kappa_1$ satisfy the left Sonine condition (\ref{eq:left_soncon_1}). Similarly, for the right sided operators, we construct the kernel with right convolution. Note that, the above example was first constructed in \cite{luchko2021general} defined on semi-infinite interval, however, such a construction verifies for the case of finite intervals too.

\vspace{0.15in}

\noindent
With regards to the function spaces, we will use the space $C_{\alpha}$ (\ref{eq:dimovski}) introduced in \cite{dimovski1966operational} and also used in \cite{luchko_finite_interval} for general fractional calculus. 


\begin{defination} \label{def:dimovski}
    For $\alpha \geq -1$ and $n \in \mathbb{N}$, the function spaces are defined as, 

    \begin{align}  \label{eq:dimovski}
        \begin{split}
            &C_\alpha^n (a, b] ~=~ \Bigl \{ f  : f^{(n)} \in   C_\alpha (a, b] \Bigr \}, \\
             &C_\alpha^n [a, b) ~=~ \Bigl \{ f  : f^{(n)} \in   C_\alpha [a, b)  \Bigr \}.
        \end{split}
    \end{align}

where,            
        \begin{equation*}
            C_\alpha (a, b] ~=~ \Bigl \{ f : (a, b] \rightarrow \mathbb{R} : f(x) = (x - a)^p f_1, p > \alpha, f_1 \in C [a, b] \Bigr \}, 
         \end{equation*}
        \begin{equation*}
            C_\alpha [a, b) ~=~ \Bigl \{ f : [a, b) \rightarrow \mathbb{R} : f(x) = (b - x)^p f_1, p > \alpha, f_1 \in C [a, b] \Bigr \}.
        \end{equation*}
\end{defination}

Note that, space $C_{-1}$ is inadequate to exclude all non-singular functions (see also \cite{arran_sonine_space}). However, we recognise that there is no need to define a function space explicitly to eliminate non-singular functions, rather satisfying the modified Sonine condition will lead to singular functions. Secondly, if there is an example of non-singular integrable function, which satisfy the modified Sonine condition then an inverse operator can be defined and results of fundamental theorems will hold, irrespective whether one considers fractional or not \cite{diethelm2020fractional}.  Thus, we will look for kernels belonging to $\mathbb{L}_n$ (Luchko class) as,

\begin{align} \label{eq:luchko_class}
    \begin{split}
         \mathbb{L}_n(a,b] = \Bigl \{ k_n, \kappa_n \in C_{-1}^n (a, b] :& \int_a^x k_n(x-t) \kappa_n (t) dt =  \frac{ (x - a)^{n-1}} { (n-1)!}, \\
         &n \in \mathbb{N}, a < x \leq b \in \mathbb{R} \Bigr \} \\
         \mathbb{L}_n[a,b) = \Bigl \{ k_n, \kappa_n \in C_{-1}^n [a, b) :& \int_x^b k_n(x-t) \kappa_n (t) dt =  \frac{ (b - x)^{n-1}} { (n-1)!}, \\
         &n \in \mathbb{N}, a \leq x < b \in \mathbb{R} \Bigr \} \\
    \end{split}
\end{align}

\vspace{0.15in}
\noindent
Now, we define the general fractional integral and derivatives.

\begin{defination}  If ($k_n, \kappa_n$) are a Sonine kernel from $\mathbb{L}_n(a,b]$, then, we define, 

 \textbf{(a)} The left-sided general fractional integral is defined with the kernel, $\kappa_n$ as, 

    \begin{equation}
        {}_a \mathcal{I}_{x}^{(\kappa_n)} f (x) := \int_{a}^{x} \kappa_n (x - s) f (s) ds
    \end{equation}

\textbf{(b)} The left-sided general fractional Riemann–Liouville derivative is defined with the kernel, $k_n$ as,

\begin{equation}
    {}_a^{RL} \mathcal{D}_{x}^{(k_n)} f (x) :=  \frac{d^n}  {d x^n} \int_{a}^{x} k_n (x - s) f (s) ds.
\end{equation}

\textbf{(c)} The left-sided general fractional Caputo derivative is defined with the kernel, $k_n$ as,

\begin{equation}
    {}_a^{C} \mathcal{D}_{x}^{(k_n)} f (x) :=   \int_{a}^{x} k_n (x - s) f^{(n)} (s)  ds.
\end{equation}

\end{defination}

\begin{defination}  If ($k_n, \kappa_n$) are a Sonine kernel from $\mathbb{L}_n[a,b)$, then, we define,

 \textbf{(a)} The right-sided general fractional integral is defined with the kernel, $\kappa_n$ as, 

    \begin{equation}
        {}_x \mathcal{I}_{b}^{(\kappa_n)} f (x) :=  \int_{x}^{b} \kappa_n (s - x) f (s) ds. 
    \end{equation} 

\textbf{(b)} The right-sided general fractional Riemann–Liouville derivative is defined with the kernel, $k_n$ as,

\begin{equation}
    {}_x^{RL} \mathcal{D}_{b}^{(k_n)} f (x) :=  (-1)^n \frac{d^n}  {d x^n} \int_{x}^{b} k_n ( s - x) f (s) ds.
\end{equation}

\textbf{(c)} The right-sided general fractional Caputo derivative is defined with the kernel, $k_n$ as,

\begin{equation}
    {}_x^{C} \mathcal{D}_{b}^{(k_n)} f (x) :=   (-1)^n \int_{x}^{b} k_n (x - s) f^{(n)} (s)  ds.
\end{equation}

\end{defination}

\vspace{0.15in}
\noindent
We state the below lemma for the relationship between the two types of derivative operators. 

\begin{lemma} \label{lm:gen_fin_rl_caputo}
    If ($k_n, \kappa_n$) are a Sonine kernel from $\mathbb{L}_n(a,b]$, $f \in C^n [a, b]$ and $x \in (a, b]$ then,

    \begin{align}
        \begin{split}
            \centering
             {}_a^{C} \mathcal{D}_{x}^{(k_n)} f (x)  ~&=~ {}_a^{RL} \mathcal{D}_{x}^{(k_n)} f (x) ~-~ \sum_{j=0}^{n-1} f^{(j)}(a) k_n^{(n-j-1)}(x-a) \\
             &=~ {}_a^{RL} \mathcal{D}_{x}^{(k_n)} \left [f (.) ~-~ \sum_{j=0}^{n-1} f^{(j)}(a) \{1\}^{j+1}_l  \right] (x) \\
        \end{split}
    \end{align}
\end{lemma}

\noindent
\textbf{Proof of lemma \ref{lm:gen_fin_rl_caputo}}

\begin{align*}
    \begin{split}
        \centering
        &{}_a^{RL} \mathcal{D}_{x}^{(k_n)} f (x) ~=~ \frac{d^n}  {d x^n} \int_{a}^{x} k_n (x - s) f (s) ds \\
        &\text{By change of variables,} \\
        &=~\frac{d^n }{ d x^n} \int_{0}^{x-a} k_n (y) f (x-y) dy \\
        &=~\frac{d^{n-1} }{ d x^{n-1}} \left [ \frac{d}  {d x}\int_{0}^{x-a} k_n (y) f (x-y) dy \right ]\\
        & \text{By Leibniz integral rule, we have,} \\
        &=~\frac{d^{n-1} }{ d x^{n-1}} \left [ \int_{a}^{x} k_n (x-s) \frac{d}  {d x} f (s) ds  ~+~ k_n(x-a) f(a) \right ]\\
        & \text{ let $g (x):= \frac{d f} {d x}$, we have,} \\
        &=~\frac{d^{n-2} }{ d x^{n-2}} \left [ \frac{d }{ d x} \int_{a}^{x} k_n (x-s) g (s) ds  ~+~ f(a) \frac{d }{ d x} k_n(x-a) \right ]\\
        & \text{By change of variables and Leibniz integral rule, we have,} \\
        &=~\frac{d^{n-2} }{ d x^{n-2}} \left [  \int_{a}^{x} k_n (x-s) \frac{d }{ d x} g (s) ds + k_n(x-a) g(a) + f(a) \frac{d }{ d x} k_n(x-a) \right ]\\
        & \text{By Induction, we have,} \\
        &=~\int_{a}^{x} k_n (x-s) \frac{d^n }{ d^n x} f (s) ds  ~+~ \sum_{j=0}^{n-1} f^{(j)}(a) k_n^{(n-j-1)}(x-a) \\
        &=~ {}_a^{C} \mathcal{D}_{x}^{(k_n)} f (x) ~+~ \sum_{j=0}^{n-1} f^{(j)}(a) k_n^{(n-j-1)}(x-a)
        \end{split}
\end{align*}

This completes the proof. Furthermore, it follows,

\begin{align*}
    \begin{split}
        \centering
        {}_a^{C} \mathcal{D}_{x}^{(k_n)} f (x) ~&=~  {}_a^{RL} \mathcal{D}_{x}^{(k_n)} f (x) ~-~ \sum_{j=0}^{n-1} f^{(j)}(a) k_n^{(n-j-1)}(x-a) \\
        &=~ {}_a^{RL} \mathcal{D}_{x}^{(k_n)} f (x) ~-~ \sum_{j=0}^{n-1} f^{(j)}(a)  {d^n \over dx^n} {}_0 \mathcal{I}_{x-a}^{(j+1)} k_n   (x) \\
        &=~ {}_a^{RL} \mathcal{D}_{x}^{(k_n)} f (x) ~-~ \sum_{j=0}^{n-1} f^{(j)}(a)  {d^n \over dx^n} \left ( k_n *  \{1\}^{j+1}_l \right )  (x) \\
         &=~ {}_a^{RL} \mathcal{D}_{x}^{(k_n)} f (x) ~-~ \sum_{j=0}^{n-1} f^{(j)}(a)  {}_a^{RL} \mathcal{D}_{x}^{(k_n)}  \{1\}^{j+1}_l  (x) \\
        &=~ {}_a^{RL} \mathcal{D}_{x}^{(k_n)} \left [f (.) ~-~ \sum_{j=0}^{n-1} f^{(j)}(a) \{1\}^{j+1}_l \right] (x)
\end{split}
\end{align*}

\vspace{0.15in}
\noindent
Similarly, for the right-sided operators, we state the below lemma for the relationship between the two types of derivative operators.

\begin{lemma} \label{lm:gen_fin_right_rl_caputo}
    If ($k_n, \kappa_n$) are a Sonine kernel from $\mathbb{L}_n[a,b)$,  $f \in C^n [a, b]$ and $x \in [a, b)$ then,

    \begin{align}
        \begin{split}
            \centering
             {}_x^{C} \mathcal{D}_{b}^{(k_n)} f (x)  ~&=~ {}_x^{RL} \mathcal{D}_{b}^{(k_n)} f (x) ~-~ \sum_{j=0}^{n-1} (-1)^{j}f^{(j)}(b) k^{(n-j-1)}(b-x) \\
             &=~ {}_x^{RL} \mathcal{D}_{b}^{(k_n)} \left [f (.) ~-~ \sum_{j=0}^{n-1} (-1)^{j-n}f^{(j)}(b) \{1\}_r^{j+1} \right] (x) \\
        \end{split}
    \end{align}
\end{lemma}

\noindent
\textbf{Proof of lemma \ref{lm:gen_fin_right_rl_caputo}}

\begin{align*}
    \begin{split}
        \centering
        &{}_x^{RL} \mathcal{D}_{b}^{(k_n)} f (x) ~=~(-1)^n {d^n \over d x^n} \int_{x}^{b} k_n (s-x) f (s) ds \\
        &\text{By change of variables,} \\
        &=~(-1)^n {d^n \over d x^n} \int_{0}^{b-x} k_n (y) f (x+y) dy \\
        &=~~(-1)^{n-1}{d^{n-1} \over d x^{n-1}} \left [ - {d \over d x}\int_{0}^{b-x} k_n (y) f (x+y) dy \right ]\\
        & \text{By Leibniz integral rule, we have,} \\
        &=~(-1)^{n-1} {d^{n-1} \over d x^{n-1}} \left [ -\int_{x}^{b} k_n (s-x) {d \over d x} f (s) ds  ~+~ k(b-x) f(b) \right ]\\
        & \text{ let $g (x):= {d f \over d x}$, we have,} \\
        &=~(-1)^{n-2}{d^{n-2} \over d x^{n-2}} \left [ -{d \over d x} \int_{x}^{b} k_n (s-x) g (s) ds  ~-~ f(b) {d \over d x} k(b-x) \right ]\\
        & \text{By change of variables and Leibniz integral rule, we have,} \\
        &=(-1)^{n-2}{d^{n-2} \over d x^{n-2}} \left [ - \int_{x}^{b} k_n (s-x) {d \over d x} g (s) ds + k(b-x) g(b) - f(b) {d \over d x} k(b-x) \right ]\\
        & \text{By Induction, we have,} \\
        &=~(-1)^{n}\int_{x}^{b} k_n (s-x) {d^n \over d^n x} f (s) ds  ~+~ \sum_{j=0}^{n-1} (-1)^{j}f^{(j)}(b) k^{(n-j-1)}(b-x) \\
        &=~ {}_x^{C} \mathcal{D}_{b}^{(k_n)} f (x) ~+~ \sum_{j=0}^{n-1} (-1)^{j} f^{(j)}(b) k^{(n-j-1)}(b-x)
        \end{split}
\end{align*}

This completes the proof. Furthermore, it follows,

\begin{align*}
    \begin{split}
        &{}_x^{C} \mathcal{D}_{b}^{(k_n)} f (x) ~=~  {}_x^{RL} \mathcal{D}_{b}^{(k_n)} f (x) ~-~ \sum_{j=0}^{n-1} (-1)^{j} f^{(j)}(b) k^{(n-j-1)}(b-x) \\
        &=~  {}_x^{RL} \mathcal{D}_{b}^{(k_n)} f (x) ~-~ \sum_{j=0}^{n-1} (-1)^{j} f^{(j)}(b) {d^n \over dx^n} {}_0 \mathcal{I}^{j+1}_{b-x} k_n (x) \\
         &=~  {}_x^{RL} \mathcal{D}_{b}^{(k_n)} f (x) ~-~ \sum_{j=0}^{n-1} (-1)^{j-n}  f^{(j)}(b) (-1)^{n} {d^n \over dx^n}  \left ( k_n * \{1\}_r^{j+1} \right)(x) \\
          &=~  {}_x^{RL} \mathcal{D}_{b}^{(k_n)} f (x) ~-~ \sum_{j=0}^{n-1}  (-1)^{j-n}  f^{(j)}(b) {}_x^{RL} \mathcal{D}_{b}^{(k_n)} \{1\}_r^{j+1} \\
        &=~ {}_x^{RL} \mathcal{D}_{b}^{(k_n)} \left [f (.) ~-~ \sum_{j=0}^{n-1} (-1)^{j-n}   f^{(j)}(b) \{1\}_r^{j+1} \right] (x)
\end{split}
\end{align*}

\vspace{0.15in}
\noindent
We state the first fundamental theorem of calculus for left-sided operators.

\begin{theorem} \label{eq:first_funda_left_operators}
    (First fundamental theorem of calculus for left-sided operators)  If pair $(k_n, \kappa_n) \in \mathbb{L}_n (a, b]$ satisfy the left modified Sonine condition for $n \in \mathbb{N}$, where $a < b \in \mathbb{R}$, then, 

    \textbf{(a)} The left-sided general fractional Riemann–Liouville derivative is defined with the kernel, $k_n$  is a left inverse of left-sided general fractional integral defined with the kernel, $\kappa_n$, then, 

    \begin{equation}
        {}_a^{RL} \mathcal{D}_{x}^{(k_n)}  {}_a \mathcal{I}_{x}^{(\kappa_n)} f (x) = f (x) ~,~ a < x \le b. 
    \end{equation} 

    \textbf{(b)} The left-sided general fractional Caputo derivative is defined with the kernel, $k_n$  is a left inverse of left-sided general fractional integral defined with the kernel, $\kappa_n$,  then, 

    \begin{equation}
       {}_a^{C} \mathcal{D}_{x}^{(k_n)}  {}_a \mathcal{I}_{x}^{(\kappa_n)} f (x) = f (x) ~,~ a < x \le b. 
    \end{equation} 

     If there exists, $\phi(x) := \left \{ \phi(x) \in L_1[a,b] : \phi(x) = {}_a \mathcal{I}_{x}^{(\kappa_n)} f (x) , x \in [a,b] \right \}$. 
\end{theorem}

\textbf{Proof of theorem \ref{eq:first_funda_left_operators}:} We split the proof into parts.

\textbf{Part (a):} First we prove the case involving Riemann–Liouville derivative. Consider,

\begin{align*}
\begin{split}
    \centering
     {}_a^{RL} \mathcal{D}_{x}^{(k_n)}  {}_a \mathcal{I}_{x}^{(\kappa_n)} f (x) ~&=~ {d^n \over d x^n} \left (  k_n * \kappa_n * f \right) (x) \\
     &=~ {d^n \over d x^n} \left ( \{1\}^n_l * f \right) (x) \\
     &=~ {d^n \over d x^n} \left ( {}_a \mathcal{I}^n_x  f \right) (x) \\
     &=~ f(x). ~~~\text{This completes the proof.}
\end{split}    
\end{align*}

\textbf{Part (b):} We  now prove the case involving Caputo derivative. Let the auxiliary function, $\phi (x) ~:=~  {}_a \mathcal{I}_{x}^{(\kappa_n)} f (x) $. Now, consider,

\begin{align*}
\begin{split}
  \centering
   {}_a^{C} \mathcal{D}_{x}^{(k_n)}  {}_a \mathcal{I}_{x}^{(\kappa_n)} f (x) ~&=~  {}_a^{C} \mathcal{D}_{x}^{(k_n)} \phi (x) \\
    &=~  {}_a^{RL} \mathcal{D}_{x}^{(k_n)} \phi (x)  ~-~ \sum_{j=0}^{n-1} \phi^{(j)}(a) k_n^{(n-j-1)}(x-a)  \\
    &=~  {}_a^{RL} \mathcal{D}_{x}^{(k_n)} {}_a \mathcal{I}_{x}^{(\kappa_n)} f (x)  \\
    &=~ f(x). 
\end{split}    
\end{align*}

Note that in the above proof we use the fact, $\phi (a) =  {}_a \mathcal{I}_{x}^{(\kappa_n)} f (a) = \displaystyle \lim_{x \to  a} \int_a^x f (x) = 0 $. Thus all subsequent derivatives are zero, i.e., $\phi^{(n)}(x) ~=~ 0$. This completes the proof.

\vspace{0.15in}
\noindent
We state the first fundamental theorem of calculus for right-sided operators.

\begin{theorem} \label{eq:first_funda_right_operators}
    (First fundamental theorem of calculus for right-sided operators)  If pair $(k_n, \kappa_n) \in \mathbb{L}_n[a,b)$ satisfy the right modified Sonine condition for $n \in \mathbb{N}$, where $a < b \in \mathbb{R}$, then, 

    \textbf{(a)} The right-sided general fractional Riemann–Liouville derivative is defined with the kernel, $k_n$  is a left inverse of right-sided general fractional integral defined with the kernel, $\kappa_n$, then, 

    \begin{equation}
        {}_{x}^{RL} \mathcal{D}_{b}^{(k_n)}  {}_x \mathcal{I}_{b}^{(\kappa_n)} f (x) = (-1)^n f (x) ~,~ a \le x < b. 
    \end{equation} 

    \textbf{(b)} The right-sided general fractional Caputo derivative is defined with the kernel, $k_n$  is a left inverse of right-sided general fractional integral defined with the kernel, $\kappa_n$, then, 

    \begin{equation}
       {}_{x}^{C} \mathcal{D}_{b}^{(k_n)}  {}_x \mathcal{I}_{b}^{(\kappa_n)} f (x) = (-1)^n f (x) ~,~ a \le x < b. 
    \end{equation} 
    
     If there exists, $\phi(x) := \left \{ \phi(x) \in L_1[a,b] : \phi(x) = {}_x \mathcal{I}_{b}^{(\kappa_n)} f (x) , x \in [a,b] \right \}$. 
\end{theorem}

\textbf{Proof of theorem \ref{eq:first_funda_right_operators}:} We split the proof into parts.

\textbf{Part (a):} First we prove the case involving Riemann–Liouville derivative. Consider,

\begin{align*}
\begin{split}
    \centering
     {}_{x}^{RL} \mathcal{D}_{b}^{(k_n)}  {}_{x} \mathcal{I}_{b}^{(\kappa_n)} f (x) ~&=~  (-1)^n{d^n \over d x^n} \left (  \int_b^x k_n (t-x) \int_b^t \kappa_n (\tau - x) f (\tau) d \tau dt \right)  \\
     &=~  (-1)^n{d^n \over d x^n} \left (  \int_x^b k_n (t-x) \int_x^b \kappa_n (\tau - x) f (\tau) d \tau dt \right)  \\
    &=~ (-1)^{n}{d^n \over d x^n} \left ( \{-1\}^n_r * f \right) (x) \\
     &=~ (-1)^{2n} {d^n \over d x^n} \left ( {}_x\mathcal{I}^n_b  f \right) (x) \\
     &=~  (-1)^{3n}f(x) ~=~ (-1)^{2n} \left({1 \over -1}\right)^n f (x) ~=~(-1)^{n} f (x). 
\end{split}    
\end{align*}

\textbf{Part (b):} We  now prove the case involving Caputo derivative. Let the auxiliary function, $\phi (x) ~:=~  {}_x \mathcal{I}_{b}^{(\kappa_n)} f (x) $. Now, consider,

\begin{align*}
\begin{split}
  \centering
   {}_x^{C} \mathcal{D}_{b}^{(k_n)}  {}_x \mathcal{I}_{b}^{(\kappa_n)} f (x) ~&=~  {}_x^{C} \mathcal{D}_{b}^{(k_n)} \phi (x) \\
    &=~  {}_x^{RL} \mathcal{D}_{b}^{(k_n)} \phi (x)  ~-~ \sum_{j=0}^{n-1} (-1)^{j} \phi^{(j)}(b) k_n^{(n-j-1)}(b-x)  \\
    &=~  {}_x^{RL} \mathcal{D}_{b}^{(k_n)} {}_x \mathcal{I}_{b}^{(\kappa_n)} f (x)  \\
    &=~ (-1)^{n} f(x). 
\end{split}    
\end{align*}

Note that in the above proof we use the fact, $\phi (b) =  {}_x \mathcal{I}_{b}^{(\kappa_n)} f (b) = \displaystyle \lim_{x \to  b} \int_x^b f (x) = 0 $. Thus all subsequent derivatives are zero, i.e., $\phi^{(n)}(x) = 0$. This completes the proof.

\vspace{0.15in}
\noindent
We state the second fundamental theorem of calculus for left-sided operators.

\begin{theorem} \label{eq:second_funda_left_operators}
    (Second fundamental theorem of calculus for left-sided operators)  If pair $(k_n, \kappa_n) \in \mathbb{L}_n(a,b]$ satisfy the left modified Sonine condition for $n \in \mathbb{N}$, where $a < b \in \mathbb{R}$, then, 

    \textbf{(a)} If there exists,  $\phi(x) := \left \{ \phi(x) \in L_1[a,b] : f(x) = {}_a \mathcal{I}_{x}^{(\kappa_n)} \phi (x) , x \in [a,b] \right \}$ 

    \begin{equation*}
         {}_a \mathcal{I}_{x}^{(\kappa_n)} {}_a^{RL} \mathcal{D}_{x}^{(k_n)}  f (x) = f (x) ~,~ a < x \le b. 
    \end{equation*} 

    \textbf{(b)} For a function, $f \in C^{n}[a, b]$, we have,

    \begin{equation*}
       {}_a \mathcal{I}_{x}^{(\kappa_n)} {}_a^{C} \mathcal{D}_{x}^{(k_n)}  f (x) = f (x) ~-~ \sum_{j=0}^{n-1} f^{(j)}(a) {\{1\}_l^{j+1}} (x) ~,~ a < x \le b. 
    \end{equation*} 
    
\end{theorem}

\textbf{Proof of theorem \ref{eq:second_funda_left_operators}:} We split the proof into parts.

\textbf{Part (a):} First we prove the case involving Riemann–Liouville derivative. Let $f (x) := {}_a \mathcal{I}_{x}^{(\kappa_n)} \phi (x)$. Consider,

\begin{align*}
\begin{split}
    \centering
     {}_a \mathcal{I}_{x}^{(\kappa_n)} {}_a^{RL} \mathcal{D}_{x}^{(k_n)}   f (x) ~&=~ {}_a \mathcal{I}_{x}^{(\kappa_n)} {}_a^{RL} \mathcal{D}_{x}^{(k_n)} {}_a \mathcal{I}_{x}^{(\kappa_n)} \phi (x) \\
     &=~ {}_a \mathcal{I}_{x}^{(\kappa_n)} \phi (x) \\
     &=~ f(x). ~~~\text{This completes the proof.}
\end{split}    
\end{align*}

\textbf{Part (b):} We now prove the case involving Caputo derivative. Consider,

\begin{align*}
\begin{split}
  \centering
    {}_a \mathcal{I}_{x}^{(\kappa_n)} {}_a^{C} \mathcal{D}_{x}^{(k_n)}  f (x) ~&=~ {}_a \mathcal{I}_{x}^{(\kappa_n)} \left [ {}_a^{RL} \mathcal{D}_{x}^{(k_n)}f(.)  ~-~  \sum_{j=0}^{n-1} f^{(j)}(a) {}_a^{RL} \mathcal{D}_{x}^{(k_{n})} {\{1\}_l^{j+1}} \right ] (x)\\
    &=~ f(x)  ~-~  \sum_{j=0}^{n-1} f^{(j)}(a) {}_a \mathcal{I}_{x}^{(\kappa_n)} {}_a^{RL} \mathcal{D}_{x}^{(k_{n})} {\{1\}_l^{j+1}} (x)\\
    &=~ f(x)  ~-~  \sum_{j=0}^{n-1} f^{(j)}(a){\{1\}_l^{j+1}} (x)
\end{split}    
\end{align*}

\vspace{0.15in}
\noindent
We state the second fundamental theorem of calculus for right-sided operators.

\begin{theorem} \label{eq:second_funda_right_operators}
    (Second fundamental theorem of calculus for right-sided operators)  If pair $(k_n, \kappa_n) \in \mathbb{L}_n[a,b)$ satisfy the right modified Sonine condition for $n \in \mathbb{N}$, where $a < b \in \mathbb{R}$, then, 

    \textbf{(a)}  If there exists, $\phi(x) := \left \{ \phi(x) \in L_1[a,b] : f(x) = {}_x \mathcal{I}_{b}^{(\kappa_n)} \phi (x) , x \in [a,b] \right \}$ 

    \begin{equation}
         {}_x \mathcal{I}_{b}^{(\kappa_n)} {}_x^{RL} \mathcal{D}_{b}^{(k_n)}  f (x) = (-1)^n f (x) ~,~ a \le x < b. 
    \end{equation} 

    \textbf{(b)}  For a function, $f \in C^{n}[a, b]$, we have,

    \begin{equation}
       {}_x \mathcal{I}_{b}^{(\kappa_n)} {}_x^{C} \mathcal{D}_{b}^{(k_n)}  f (x) =  (-1)^{n} f (x) - \sum_{j=0}^{n-1} (-1)^{j} f^{(j)}(b) {\{1\}_r^{j+1}} (x) ,~ a \le x < b. 
    \end{equation} 
    
\end{theorem}

\textbf{Proof of theorem \ref{eq:second_funda_right_operators}:} We split the proof into parts.

\textbf{Part (a):} First we prove the case involving Riemann–Liouville derivative. Let $f (x) := {}_x \mathcal{I}_{b}^{(\kappa_n)} \phi (x)$. Consider,

\begin{align*}
\begin{split}
    \centering
     {}_x \mathcal{I}_{b}^{(\kappa_n)} {}_x^{RL} \mathcal{D}_{b}^{(k_n)}   f (x) ~&=~ {}_x \mathcal{I}_{b}^{(\kappa_n)} {}_x^{RL} \mathcal{D}_{b}^{(k_n)} {}_x \mathcal{I}_{b}^{(\kappa_n)} \phi (x) \\
     &=~ {}_x \mathcal{I}_{b}^{(\kappa_n)} (-1)^n \phi (x) \\
     &=~ (-1)^n f(x). ~~~\text{This completes the proof.}
\end{split}    
\end{align*}

\textbf{Part (b):} We now prove the case involving Caputo derivative. Consider,

\begin{align*}
\begin{split}
  \centering
    {}_x \mathcal{I}_{b}^{(\kappa_n)} {}_x^{C} \mathcal{D}_{b}^{(k_n)}  f (x) &= {}_x \mathcal{I}_{b}^{(\kappa_n)} \left [ {}_x^{RL} \mathcal{D}_{b}^{(k_n)}f(.)  -  \sum_{j=0}^{n-1} (-1)^{j-n} f^{(j)}(b) {}_x^{RL} \mathcal{D}_{b}^{(k_{n})} {\{1\}_r^{j+1}} \right ] (x)\\
    &=~ (-1)^n f(x)  ~-~  \sum_{j=0}^{n-1}  (-1)^{j-n}  f^{(j)}(b) {}_x \mathcal{I}_{b}^{(\kappa_n)} {}_x^{RL} \mathcal{D}_{b}^{(k_{n})} {\{1\}_r^{j+1}} (x)\\
    &=~ (-1)^n f(x)  ~-~  \sum_{j=0}^{n-1}  (-1)^{j} f^{(j)}(b){\{1\}_r^{j+1}} (x)
\end{split}    
\end{align*}

\vspace{0.15in}
\noindent
Integration by parts is a key tool for both mathematical analysis and numerical methods, in this view we extend it for general fractional operators as below.

\begin{theorem} \label{th:inte_part_rl_n}
    General fractional integration by parts for Riemann–Liouville type derivative for a general $n \in \mathbb{N}$.

\begin{align}
    \begin{split}
         \int_{a}^{b} f(x)  \left ({}_a^{RL} \mathcal{D}_{x}^{(k_n)} y(x) \right) dx ~  =~  \int_{a}^{b} \left ( {}_x^{RL} \mathcal{D}_{b}^{(k_n)} f(x) \right) y(x) dx
    \end{split}
\end{align}

\end{theorem}

\noindent
\textbf{Proof of theorem \ref{th:inte_part_rl_n}:} Consider, 

\begin{align*}
    \begin{split}
         \int_{a}^{b}& f(x)  \left ({}_a^{RL} \mathcal{D}_{x}^{(k_n)} y(x) \right) dx ~=~  \int_{a}^{b} f(x) {d^n \over dx^{n} } \int_{a}^{x} k_{n} (x-s) y(s) ds dx
    \end{split}
\end{align*}

By integration by parts,

\begin{align*}
    \begin{split}
          &=~  \left [ f(x) {d^{n-1} \over dx^{n-1} } \int_{a}^{x} k_{n} (x-s) y(s) ds \right ]_{x=a}^{x=b} \\
          &~~~~~~~~-\int_{a}^{b} f^{(1)}(x) {d^{n-1} \over dx^{n-1} } \int_{a}^{x} k_{n} (x-s) y(s) ds dx \\
        &=~   f(b) {d^{n-1} \over dx^{n-1} } \int_{a}^{b} k_{n} (b-s) y(s) ds  \\
        &~~~~~~~~~-~\int_{a}^{b} f^{(1)}(x) {d^{n-1} \over dx^{n-1} } \int_{a}^{x} k_{n} (x-s) y(s) ds dx 
    \end{split}
\end{align*}

By repeated integration by parts,

\begin{align*}
    \begin{split}
        &=~  \sum_{j=0}^{n-1} (-1)^{j}  f^{(j)} (b) {d^{n-j-1} \over dx^{n-j-1} } \int_{a}^{b} k_{n} (b-s) y(s) ds \\
        &~~~~~~~~+~(-1)^n\int_{a}^{b} f^{(n)}(x) \int_{a}^{x} k_{n} (x-s) y(s) ds dx \\
        &=~ \sum_{j=0}^{n-1} (-1)^{j}  f^{(j)} (b) {d^{n-j-2} \over dx^{n-j-2} } \left ( {d \over dx } \int_{a}^{b} k_{n} (b-s) y(s) ds \right ) \\
        &~~~~~~~~+~(-1)^n\int_{a}^{b} f^{(n)}(x) \int_{a}^{x} k_{n} (x-s) y(s) ds dx
    \end{split}
\end{align*}

By Leibniz integral rule,

\begin{align*}
    \begin{split}
        &=~  \sum_{j=0}^{n-1} (-1)^{j} f^{(j)} (b) {d^{n-j-2} \over dx^{n-j-2} } \left (  \int_{a}^{b} k_{n}^{(1)} (b-s) y(s) ds + k_n(b-a) y(a) \right ) \\ 
        &~~~~~~~~+~(-1)^n\int_{a}^{b} f^{(n)}(x) \int_{a}^{x} k_{n} (x-s) y(s) ds dx \\
         &=~  \sum_{j=0}^{n-1} (-1)^{j} f^{(j)} (b) {d^{n-j-3} \over dx^{n-j-3} } \left ( {d \over dx } \int_{a}^{b} k_{n}^{(1)} (b-s) y(s) ds  \right )\\ 
         &~~~~~~~~~+~(-1)^n\int_{a}^{b} f^{(n)}(x) \int_{a}^{x} k_{n} (x-s) y(s) ds dx
    \end{split}
\end{align*}

Repeated application of Leibniz integral rule,

\begin{align*}
    \begin{split}
        &=~  \sum_{j=0}^{n-1} (-1)^{j}f^{(j)} (b) \int_{a}^{b} k_{n}^{(n-j-1)} (b-s) y(s) ds   \\
        &~~~~~~~~~+~(-1)^n\int_{a}^{b} f^{(n)}(x) \int_{a}^{x} k_{n} (x-s) y(s) ds dx \\
        &=~   \int_{a}^{b} \sum_{j=0}^{n-1} (-1)^{j} f^{(j)} (b) k_{n}^{(n-j-1)} (b-s) y(s) ds  \\
        &~~~~~~~~~+~(-1)^n\int_{a}^{b} f^{(n)}(x) \int_{a}^{x} k_{n} (x-s) y(s) ds dx 
    \end{split}
\end{align*}

Change of order of integration,

\begin{align*}
    \begin{split}
        &=~   \int_{a}^{b} \sum_{j=0}^{n-1} (-1)^{j} f^{(j)} (b) k_{n}^{(n-j-1)} (b-s) y(s) ds   \\
        &~~~~~~~~~+~(-1)^n\int_{a}^{b} y(s) \int_{s}^{b} k_{n} (x-s) f^{(n)}(x) dx ds \\
        &=   \int_{a}^{b} y(s)\left [\sum_{j=0}^{n-1} (-1)^{j} f^{(j)} (b) k_{n}^{(n-j-1)} (b-s)     +(-1)^n\int_{s}^{b} k_{n} (x-s) f^{(n)}(x) dx \right ]ds \\
        &=~  \int_{a}^{b} \left ( {}_x^{RL} \mathcal{D}_{b}^{(k_n)} f(x) \right) y(x) dx
    \end{split}
\end{align*}


\section{Basis function} \label{sec:basis}

Jacobi polynomials forms a basis function and satisfies orthogonality with respect to weighted inner product. For more details over orthogonal polynomials, refer \cite{chihara2011introduction}. We denote, $P_n^{\alpha, \beta} (x)$ as the Jacobi polynomial (\ref{eq:jac})\cite{hesthaven2007spectral}. Note, that for $\alpha = \beta = 0$ is the Legendre polynomial.

\begin{equation} \label{eq:jac}
    P_n^{\alpha, \beta} (x) ~=~ \frac{1}{2^n} \sum_{k=0}^{n} \binom{n+\alpha}{k}  \binom{n+\beta}{n-k} (x-1)^{n-k} (x+1)^{k}, ~ \alpha, \beta > -1, x \in [-1, 1].
\end{equation}

\noindent
The Jacobi polynomial follows the symmetric relationship as (\ref{eq:jac_sym}) \cite{hesthaven2007spectral}.

\begin{equation} \label{eq:jac_sym}
    P_n^{\alpha, \beta} (x) ~=~ (-1)^n P_n^{\beta, \alpha} (-x)
\end{equation}

\noindent
For the case of Jacobi polynomial, we have the orthogonality relation as,

\begin{equation}
    \left ( P_n^{\alpha, \beta},  P_m^{\alpha, \beta}   \right)_{L^2_w [-1, 1]} ~=~ \frac{2^{\alpha+\beta+1} }{(2n + \alpha + \beta + 1)n!} \frac{\Gamma(n+\alpha+1)\Gamma(n+\beta+1) }{\Gamma(n+\alpha+\beta+1)} \delta_{mn}
\end{equation}

where, $\delta_{mn}$ denotes the Kronecker delta. Indeed, we shall denote the orthogonality constant, $\gamma_n = \|  P_n^{\alpha, \beta} \|^2_{L^2_w [-1, 1]} =  \left ( P_n^{\alpha, \beta},  P_n^{\alpha, \beta}   \right)_{L^2_w [-1, 1]} $, unless otherwise stated explicitly.

\noindent
The Jacobi Polynomials satisfy the following three term recurrence (\ref{eq:jac_three}) \cite{hesthaven2007spectral} for $n \geq 0$,

\begin{align} \label{eq:jac_three}
    \begin{split}
        x  P_n^{\alpha, \beta} (x) &= a_{n-1, n}^{\alpha, \beta} P_{n-1}^{\alpha, \beta} (x) + a_{n, n}^{\alpha, \beta} P_{n}^{\alpha, \beta} (x) + a_{n+1, n}^{\alpha, \beta} P_{n+1}^{\alpha, \beta} (x) \\
        \text{where,} \\
        &  a_{n-1, n}^{\alpha, \beta} = \frac{2(n+\alpha)(n+\beta)}{(2n + \alpha + \beta +1)(2n + \alpha + \beta)} \\
        &  a_{n, n}^{\alpha, \beta} = \frac{\alpha^2 - \beta^2}{(2n + \alpha + \beta +2)(2n + \alpha + \beta)} \\
        &  a_{n+1, n}^{\alpha, \beta} = \frac{2(n+1)(n+\alpha+\beta+1)}{(2n + \alpha + \beta +2)(2n + \alpha + \beta+1)} 
    \end{split}
\end{align}

where, for $n = 0, a_{-1, n}^{\alpha, \beta} = 0$ and to start the three term recurrence,  we have,

\begin{align}
    \begin{split}
        & P_0^{\alpha, \beta} (x) = 0, \\
         & P_1^{\alpha, \beta} (x) = \frac{1}{2} \left( \alpha + \beta + 2 \right) x ~+~  \frac{1}{2} \left( \alpha - \beta \right)
    \end{split}
\end{align}

\noindent
We denote, $\Tilde{P}_n^{\alpha, \beta} (x)$ as the shifted Jacobi polynomial (\ref{eq:shit_jac}) for $x \in [0, 1]$ obtained via an affine transformation. 

\begin{equation} \label{eq:shit_jac}
    \Tilde{P}_n^{\alpha, \beta} (x) ~=~  P_n^{\alpha, \beta} (2x - 1)
\end{equation}

We shall denote the orthogonality constant as $\Tilde{\gamma}_n = \|  \Tilde{P}_n^{\alpha, \beta} \|^2_{L^2_w [0, 1]} =  \left ( \Tilde{P}_n^{\alpha, \beta},  \Tilde{P}_n^{\alpha, \beta}   \right)_{L^2_w [0, 1]} $, for the case of shifted Jacobi polynomials. 

Owing to the above properties, Jacobi polynomials and its special cases such as Legendre and Chebeshev polynomials has been a popular choice in construction of spectral methods \cite{hesthaven2007spectral, gottlieb1977numerical, canuto1988spectral}.

\subsection{Jacobi convolution series}

One of the overwhelming issue on using Jacobi polynomials for fractional differential equations (also general fractional differential equations) which are convolution type operators that, we have to compute convolution resulting in a full matrix, thereby limiting accuracy and computational efficiency of the method. A better approach is to construct basis functions, such that the fractional (or general) derivative is a power series. Such functions often have non-polynomial structure. 

In the realm of fractional derivatives, in  \cite{zayernouri2013fractional}, Jacobi Poly-fractonomials were introduced (see also \cite{zayernouri2015tempered}) and in \cite{chen2016generalized}, generalised Jacobi functions were introduced. Both of these functions have non-polynomial behaviour. Furthermore, the fractional derivative of either functions is a power series, if one uses a suitable test function, then orthogonality holds (with respect to weighted inner product) and we get a diagonal stiffness matrix. This not only results in an accurate method but also an efficient scheme. 

The central idea for construction of either Jacobi Poly-fractonomials or generalised Jacobi functions relies on a key result by \citet{askey1969integral} (\ref{eq:askey}). Notice, that as a result of  convolution, on the left hand side of (\ref{eq:askey}) we have a Jacobi polynomial. Indeed, this fact used in construction of Jacobi Poly-fractonomials \cite{zayernouri2013fractional} or generalised Jacobi functions \cite{chen2016generalized}.

    \begin{equation} \label{eq:askey}
        (1+x)^{\beta + \mu} \frac{P_n^{\alpha-\mu, \beta+\mu}(x)}{P_n^{\alpha-\mu, \beta+\mu}(-1)} = \frac{\Gamma(\beta+\mu +1)}{\Gamma(\beta+1)\Gamma(\mu)} \int_{-1}^x (1 +y)^{\beta} \frac{P_n^{\alpha, \beta}(y)}{P_n^{\alpha, \beta}(-1)} (x-y)^{\mu -1} dy
    \end{equation}

\noindent
We would like to further extend this idea of seeking for functions such that the general fractional derivative of such function is a power series, explicitly written. Although, the result (\ref{eq:askey}) is powerful, but in case of general fractional operators the kernel is arbitrary, hence, the difficulty in applying (\ref{eq:askey}) for all such kernel belonging to $\mathbb{L}_n$ (and not just specific examples). Thus, to make it a general method for any choice of kernel belonging to $\mathbb{L}_n$ we introduce the Jacobi convolution series as basis functions and in the next section we construct an efficient Petrov-Galerkin scheme, where we obtain a diagonal stiffness matrix. Indeed, our method, using Jacobi convolution series is applicable to fractional derivatives too, since they are a special case of general fractional derivative. In this view, Jacobi convolution series can be regarded as generalisation of generalised Jacobi functions  \cite{chen2016generalized}, while in \cite{chen2016generalized} is shown that Jacobi Poly-fractonomials \cite{zayernouri2013fractional} is a special case of generalised Jacobi functions. Needless to mention, generalised Jacobi functions for a suitable choice of parameter leads to Jacobi polynomials.

We now introduce Jacobi convolution series and subsequently prove form a basis function.

\begin{defination} \label{def:jconv} (Left Jacobi convolution series)
    We define left Jacobi convolution series ($\phi_n(x)$) as,
    \begin{equation} \label{eq:jconv}
        \phi_n(x) := \int_0^x \kappa(x-t) \Tilde{P}_n^{\alpha, \beta} (t) dt  , ~x \in [0, 1],~ \alpha, \beta > -1, \forall n \in \mathbb{N} \cup \{0\}
    \end{equation}
    where $\kappa \in \mathbb{L}_m (0,1]$ satisfies the left modified Sonine condition and $\Tilde{P}_n^{\alpha, \beta} (x)$ is the shifted Jacobi polynomial. 
\end{defination}

By virtue of our construction, we have, 

\begin{equation}
     \phi_n(0) = \lim_{x\to 0} \int_0^x \kappa(x-t) \Tilde{P}_n^{\alpha, \beta} (t) dt = 0, ~ \forall n \in \mathbb{N} \cup \{0\}
\end{equation}

Our construction was motivated by the fact, the left-sided general fractional derivative of left Jacobi convolution series is a shifted Jacobi polynomial, shown as, 

\begin{equation} \label{eq:gfd_jconv}
    {}_0^{RL} D_x^k \phi_n ~=~ \frac{d}{dx} k * \kappa * \Tilde{P}_n^{\alpha, \beta} (x) ~=~\frac{d}{dx} \{1\} * \Tilde{P}_n^{\alpha, \beta} (x)  ~=~   \Tilde{P}_n^{\alpha, \beta} (x),  x \in (0, 1]
\end{equation}

\noindent
Similarly, we define the right Jacobi convolution series as (\ref{eq:jrconv}). From the context, it should be clear, if $\phi_n$ denotes the left or right Jacobi convolution series.

\begin{defination} \label{def:jrconv} (Right Jacobi convolution series)    We define right Jacobi convolution series ($\phi_n(x)$) as,
    \begin{equation} \label{eq:jrconv}
        \phi_n(x) := \int_x^1 \kappa(t-x) \Tilde{P}_n^{\alpha, \beta} (t) dt  , ~x \in [0, 1],~ \alpha, \beta > -1,  \forall n \in \mathbb{N} \cup \{0\}
    \end{equation}
    where $\kappa \in \mathbb{L}_m [0,1)$ satisfies the right modified Sonine condition and $\Tilde{P}_n^{\alpha, \beta} (x)$ is the shifted Jacobi polynomial.
\end{defination}

Again, the right Jacobi convolution series is zero at $x=1$, since, 

\begin{equation}
    \phi_n(1) = \lim_{x\to 1} \int_x^1 \kappa(t-x) \Tilde{P}_n^{\alpha, \beta} (t) dt = 0, \forall n \in \mathbb{N} \cup \{0\}
\end{equation}

Linear independence of $\phi_n$ (for either case) can be shown trivially, by considering $\sum_{i=0}^n c_i \phi_i = \kappa * \sum_{i=0}^n c_i \Tilde{P}_i^{\alpha, \beta} (x) = 0$. Since $\kappa$ is arbitrary, indeed it has the only solution of $\{c_i\}_{i=0}^{n} = 0$.

\begin{remark}
It is to be noted that, the Jacobi convolution series are non-polynomial functions, in general. 
\end{remark}
    
The other examples of non-polynomial functions prevalent in fractional calculus literature are the Jacobi poly-fractonomials \cite{zayernouri2013fractional} \cite{zayernouri2015tempered} and generalized Jacobi functions \cite{chen2016generalized}, both of them being basis functions and used for construction of spectral methods for fractional differential equations. 

A notable mention is the development of the Taylor series and the generalized convolution Taylor formula in \cite{luchko2022convolution} for the case of general fractional operators (see also \cite{luchko2021special}).

\vspace{0.05in}\noindent
\textbf{Example:} For an example of Jacobi convolution series; we first construct the Sonine kernel (see \cite{sonine1884generalisation}) as,

\begin{align} \label{eq:son_pair}
    \begin{split}
         k(x) &= \frac{x^{-\alpha}}{\Gamma(1-\alpha)} \sum_{k=0}^N a_k x^k \\
          \kappa(x) &= \frac{x^{\alpha-1}}{\Gamma(\alpha)} \sum_{k=0}^N b_k x^k 
    \end{split}
\end{align}

where, the coefficients follow the relationship, 

\begin{align}
    \begin{split}
        &a_0 b_0 ~=~ 1,~  k = 0 \\
        &\sum_{k=1}^N \Gamma(k+1-\alpha) \Gamma(\alpha + N - k) a_{N-k} b_k = 0,~  k = \{1, 2, \dots, N \}
    \end{split}
\end{align}

\begin{figure}[!htb]
    \centering
    \includegraphics[width=0.6\linewidth]{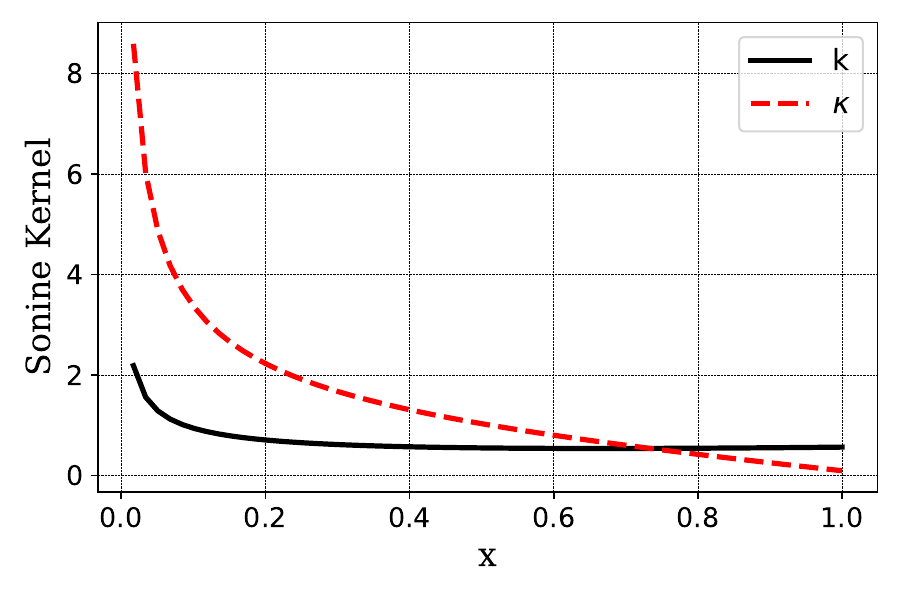}
    \caption{Sonine Kernel obtained using (\ref{eq:son_pair}) for $\alpha = 0.5$ and $a = \{0.5, 0.25, 0.25\}$ results in $b = \{2, -1, -0.83333\}$ are singular functions with singularity at $x=0$}
    \label{fig:sonine_kernel}
\end{figure}

\begin{figure} 
\centering
    \subfloat[$n=0$]{{\includegraphics[trim=0.3cm 0.3cm 0.3cm 0.3cm,clip, width=0.4\textwidth]{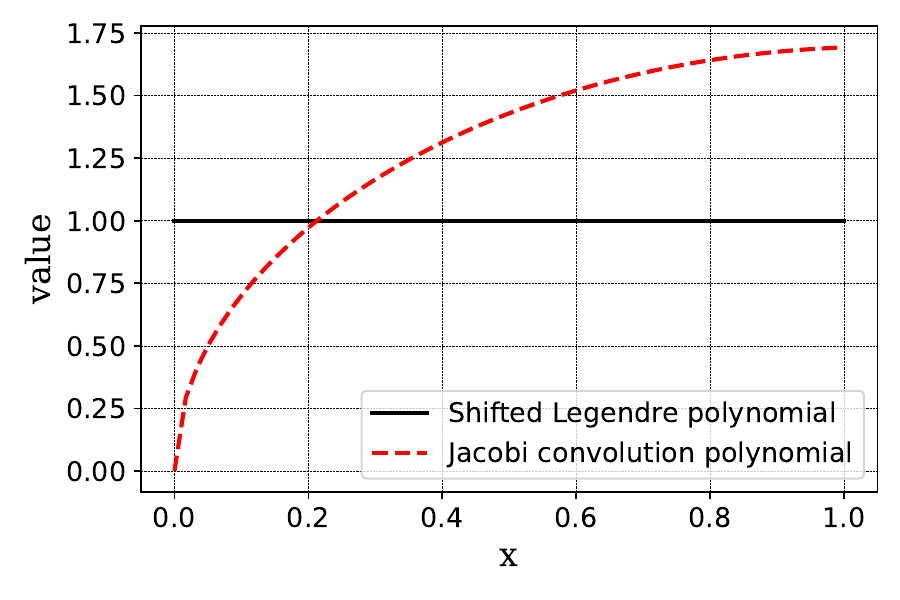}}}%
   \subfloat[$n=1$]{{\includegraphics[trim=0.3cm 0.3cm 0.3cm 0.3cm,clip, width=0.4\textwidth]{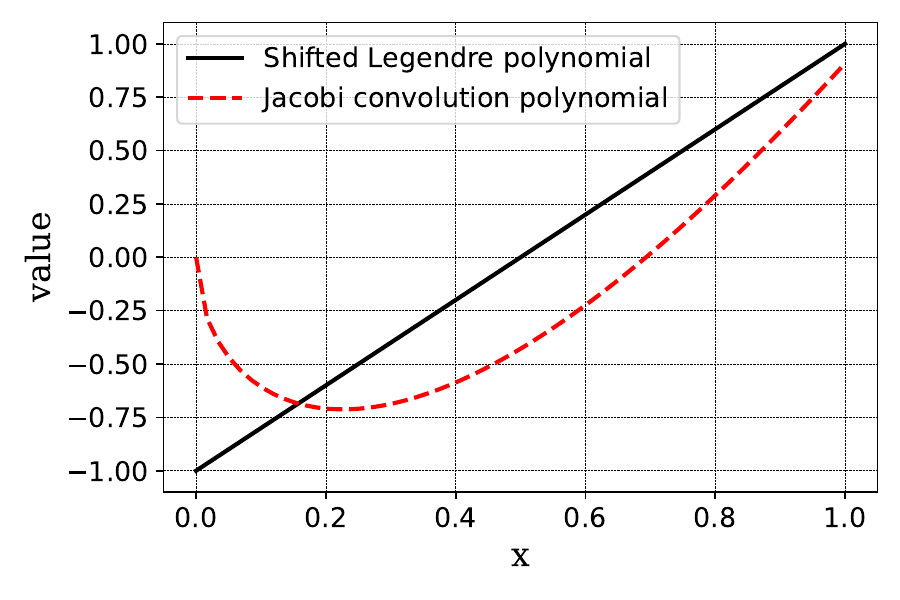}}} \\
   \subfloat[$n=2$]{{\includegraphics[trim=0.3cm 0.3cm 0.3cm 0.3cm,clip, width=0.4\textwidth]{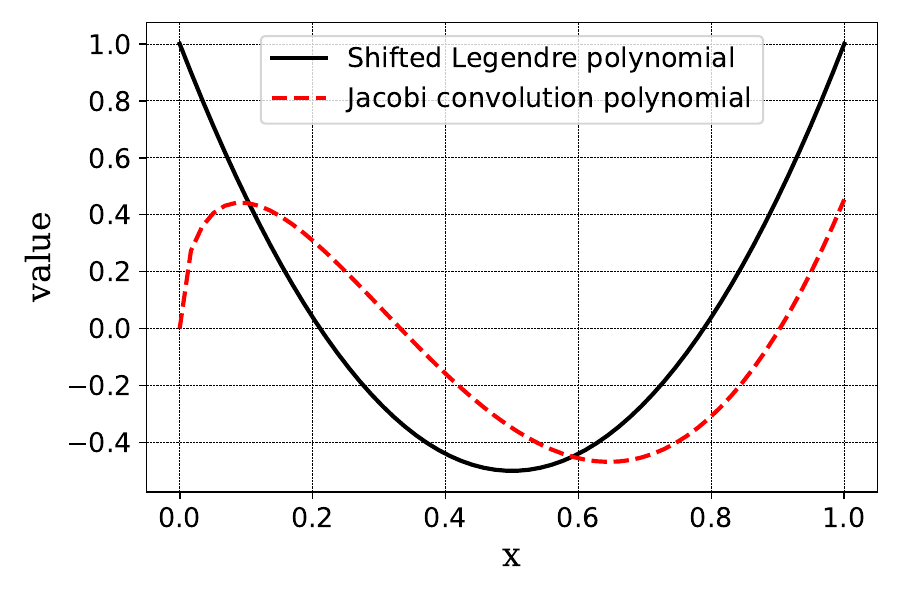}}}%
    \subfloat[$n=3$]{{\includegraphics[trim=0.3cm 0.3cm 0.3cm 0.3cm,clip, width=0.4\textwidth]{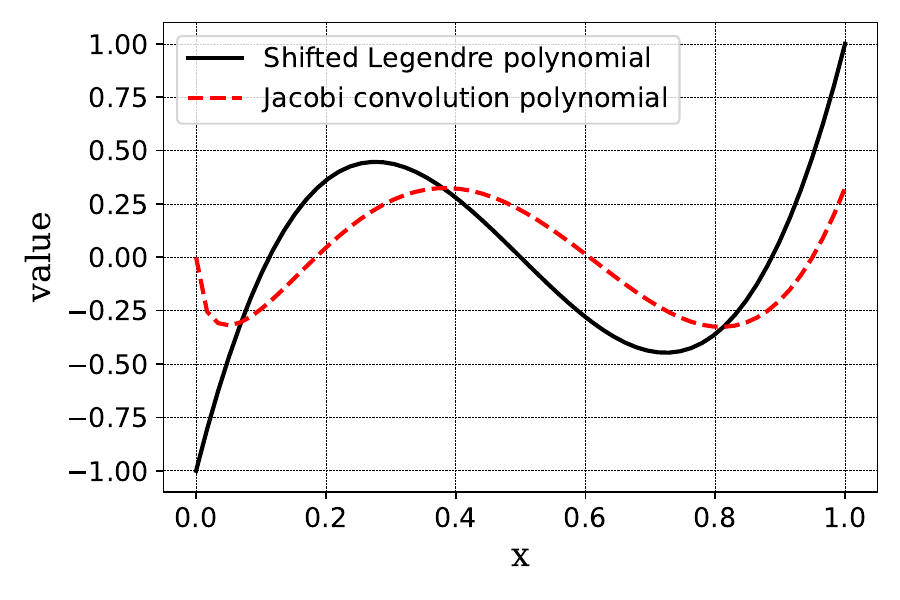}}}\\
   \subfloat[$n=4$]{{\includegraphics[trim=0.3cm 0.3cm 0.3cm 0.3cm,clip, width=0.4\textwidth]{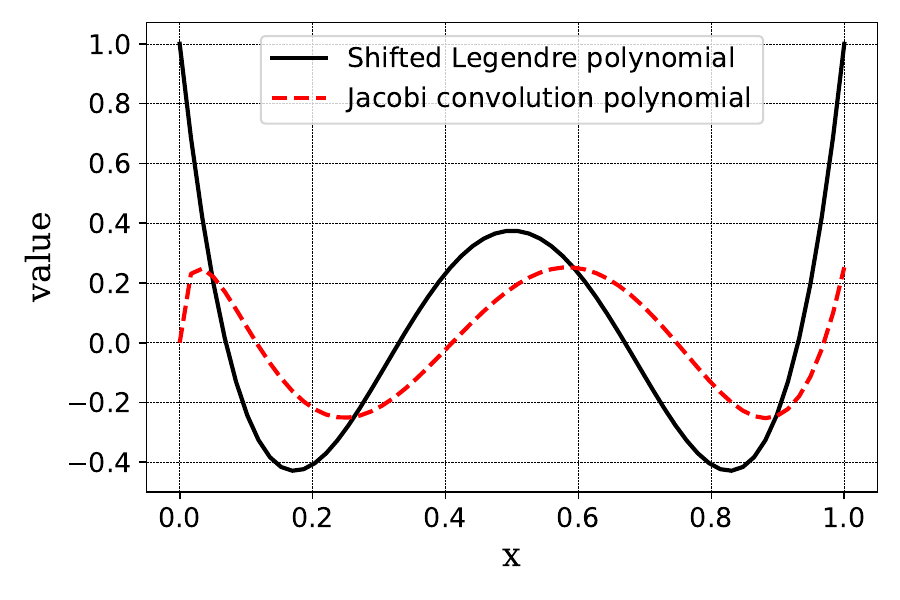}}}%
   \subfloat[$n=5$]{{\includegraphics[trim=0.3cm 0.3cm 0.3cm 0.3cm,clip, width=0.4\textwidth]{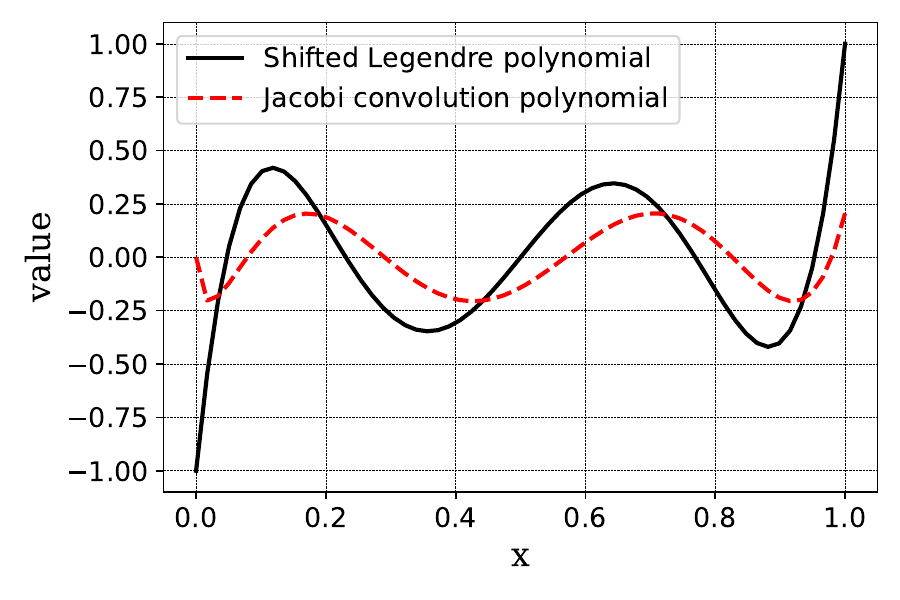}}}\\
    \subfloat[$n=6$]{{\includegraphics[trim=0.3cm 0.3cm 0.3cm 0.3cm,clip, width=0.4\textwidth]{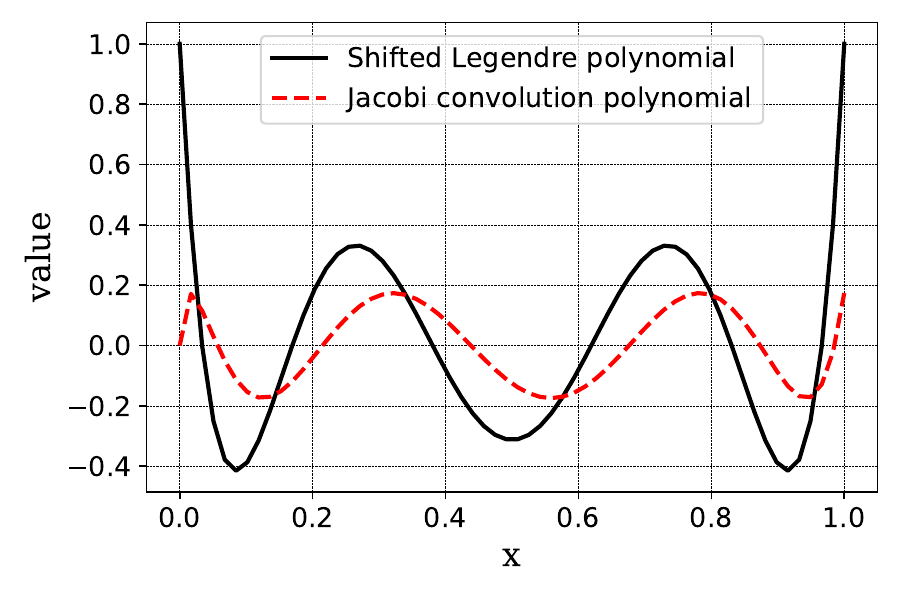}}}%
   \subfloat[$n=7$]{{\includegraphics[trim=0.3cm 0.3cm 0.3cm 0.3cm,clip, width=0.4\textwidth]{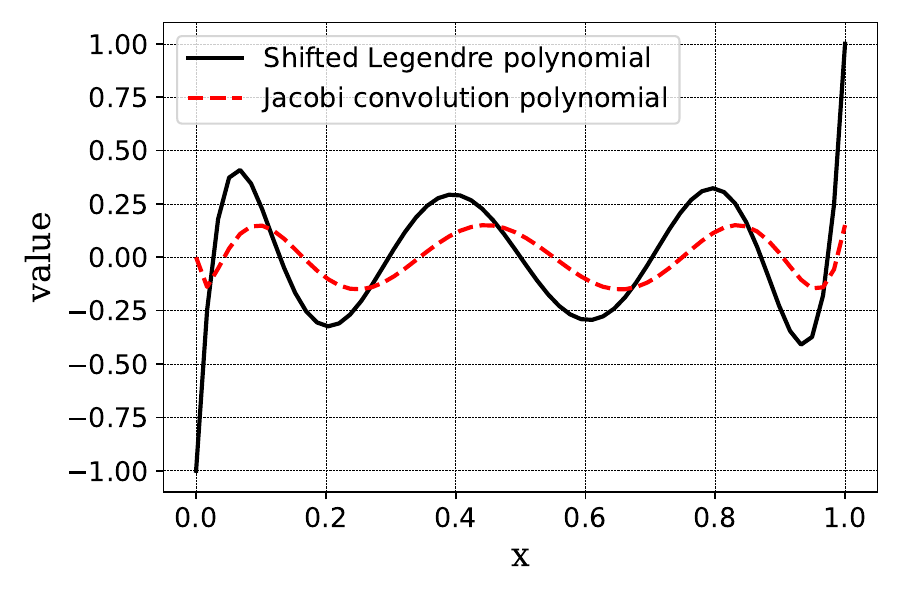}}}%
    \caption{The first eight left Jacobi convolution series (with $\alpha =\beta =0$) obtained using for the Sonine pair using (\ref{eq:son_pair}) for $\alpha = 0.5$ and $a = \{0.5, 0.25, 0.25\}$ results, $b = \{2, -1, -0.83333\}$  and the shifted Legendre polynomials}
    \label{fig:jcpl}
\end{figure}

For $\alpha = 0.5$ and $a = \{0.5, 0.25, 0.25\}$ results, $b = \{2, -1, -0.83333\}$ using the above relationship; fig.~\ref{fig:sonine_kernel} is a plot, note that the these kernel have singularity at $x=0$. Using the obtained Sonine kernel, we now plot (fig.\ref{fig:jcpl}) the left Jacobi convolution series (\ref{eq:jconv}) with $\alpha = \beta = 0$, corresponding to shifted Legendre polynomials.

\begin{theorem} \label{th:dense}
    Left Jacobi convolution series $\{\phi_n\}_{n=0}^\infty$ (\ref{eq:jconv}) form a basis in infinite dimensional Hilbert space.
    
\end{theorem}

\noindent
\textbf{Proof of theorem \ref{th:dense}:} Let, $\phi_n = \kappa *  \Tilde{P}_n^{\alpha, \beta}  = \int_0^x \kappa(x-t) \Tilde{P}_n^{\alpha, \beta} (t) dt$ and $f(x) \in L_w^2[0,1]$. 

Consider,

\begin{align}
    \begin{split}
       \left \| f - \sum_{n=0}^N a_n \phi_n \right \|_{L_w^2[0,1]} ~&=~    \left \| f - \kappa * \sum_{n=0}^N a_n  \Tilde{P}_i^{\alpha, \beta} \right \|_{L_w^2[0,1]} \\
        &=~    \left \|  \kappa * \left(g -  \sum_{n=0}^N a_n  \Tilde{P}_i^{\alpha, \beta} \right) \right \|_{L_w^2[0,1]} \\
        &\text{By Cauchy-Schwatz, we have} \\
        &\le ~  \left \| \kappa  \right \|_{L_w^2[0,1]} \left \| g -  \sum_{n=0}^N a_n  \Tilde{P}_i^{\alpha, \beta} \right \|_{L_w^2[0,1]} \\
        &\le ~\left \| g -  \sum_{n=0}^N a_n  \Tilde{P}_i^{\alpha, \beta}  \right \|_{L_w^2[0,1]} \\
    \end{split}
\end{align}

Therefore, when $N\to \infty$, 

\begin{equation}
    \lim_{N \to \infty}  \left  \| f - \sum_{n=0}^N a_n \phi_n \right \|_{L_w^2[0,1]} ~\le~  \lim_{N \to \infty} \left  \| g -  \sum_{n=0}^N a_n  \Tilde{P}_i^{\alpha, \beta} \right \|_{L_w^2[0,1]} \to 0, 
\end{equation}

by Weierstrass's theorem.  This completes the proof.

\begin{theorem} \label{th:dense2}
    Right Jacobi convolution series $\{\phi_n\}_{n=0}^\infty$ (\ref{eq:jrconv}) forms a basis in infinite dimensional Hilbert space.
\end{theorem}

\noindent
\textbf{Proof of theorem \ref{th:dense2}:} The proof is omitted, since it follows the same ideas as theorem \ref{th:dense}.

\section{Petrov-Galerkin scheme for general fractional derivative} \label{sec:pg}

For an efficient construction of a numerical scheme, we use (\ref{eq:gfd_jconv}) which leads to a Petrov-Galerkin scheme. We illustrate our construction for the boundary value problem (\ref{eq:bc_1}). 

\begin{align} \label{eq:bc_1}
\begin{split}
     &{}_0^{RL} D_x^{(k)} f (x) ~=~ g (x), ~ x \in (0, 1), \\
     &f(0) = 0 , ~ f(1) = b.
\end{split}
\end{align}

Since $\phi_n$ defined in (\ref{eq:jconv}) is a basis function, we construct a space as, 

\begin{equation}
    U :=  \text{span} ~ \Bigl \{ \phi_n : \phi_n(0) = 0,  n \in \mathbb{N} \cup \{0\}  \Bigr \}
\end{equation}

We construct the space of test function, $V$ as,

\begin{equation}
    V :=  \text{span} ~ \Bigl \{   \Tilde{P}_n^{\alpha, \beta} , \alpha, \beta > -1, n \in \mathbb{N} \cup \{0\}  \Bigr \} ,
\end{equation}

where, $ \Tilde{P}_n^{\alpha, \beta}$ is a shifted Jacobi polynomial. As a result, we obatin a bilinear form of (\ref{eq:bc_1}), for $f \in U$ and $v \in V$ as, 

\begin{equation}
    a (f, v) ~:=~ \left( {}_0^{RL} D_x^{(k)} f  ~,~   v  \right)_{L^2_w [0,1]} ~=~  \left( g ~,~  v  \right)_{L^2_w [0,1]} 
\end{equation}

\noindent
For the numerical approximation of $f$, we seek the solution ($f_N$) of the form,

\begin{equation} \label{eq:jac_galerkin}
    f_N (x) ~=~ \sum_{n=0}^N \hat{f}_n \phi_n (x) ~=~ \sum_{n=0}^N  \hat{f}_n  \left( \kappa * \Tilde{P}_n^{\alpha, \beta} \right) (x)
\end{equation}

where, $f_N \in U_N$ and $U_N \subset U$ is a finite dimensional sub-space, dense in $U$ and $ \hat{f}_n$ are the expansion coefficients. Furthermore, $V_N \subset V$   is also finite dimensional sub-space, dense in $V$. Thus, we seek the numerical approximation, $f_N \in U_N$ and $v_N \in V_N$, such that,

\begin{equation} \label{eq:dsc_biline}
    a (f_N, v_N) ~:=~ \left( {}_0^{RL} D_x^{(k)} f_N  ~,~   v_N  \right)_{L^2_w [0,1]} ~=~  \left( g ~,~  v_N  \right)_{L^2_w [0,1]} 
\end{equation}

Indeed, it enforces the the residual, $R_N$ (\ref{eq:resi}) to be $L^2$ orthogonal to every $v_N \in V_N$.  

\begin{equation} \label{eq:resi}
    R_N := {}_0^{RL} D_x^{(k)} f_N  - g
\end{equation}

Plugging (\ref{eq:jac_galerkin}) in (\ref{eq:dsc_biline}) and using (\ref{eq:gfd_jconv}), we have,

\begin{equation}
    \sum_{n=0}^N   \hat{f}_n \left( \Tilde{P}_n^{\alpha, \beta} ~,~   \Tilde{P}_m^{\alpha, \beta} \right)_{L^2_w [0,1]} ~=~  \left( g ~,~   \Tilde{P}_m^{\alpha, \beta} \right)_{L^2_w [0,1]} , ~~ \forall m \in [0, N]
\end{equation}

The weight function, $w(x)$ is selected such that orthogonality holds. Therefore, we evaluate the coefficients as, 

\begin{equation} \label{eq:PG_const}
     \hat{f}_n ~=~  \frac{1}{\Tilde{\gamma}_n} \left( g ~,~   \Tilde{P}_n^{\alpha, \beta}  \right)_{L^2_w [0,1]}  
\end{equation}

where, $\Tilde{\gamma}_{n} = \|   \Tilde{P}_n^{\alpha, \beta} \|^2_{L^2_w[0,1]} $,  is the orthogonality constant for shifted Jacobi polynomials. The boundary condition can be applied using Tau approach or Lifting.

\vspace{0.05in}\noindent
\textbf{Example:} For a numerical example of solving (\ref{eq:bc_1}), we constructed the Sonine kernel following (\ref{eq:son_pair}), where for $\alpha = 0.5$ and $a = \{0.5, 0.25, 0.25\}$ results in $b = \{2, -1, -0.83333\}$ and fig.~\ref{fig:sonine_kernel} is a plot. Although, the kernels have singularity at $x=0$; $\phi_n (0) = \displaystyle \lim_{x \to 0} \int_{0}^x \kappa(x-t) \Tilde{P}_n^{0, 0}(t) dt = 0 $.

We consider two cases: 

(a) We solve the problem (\ref{eq:ex_1}), which has the analytical solution as $f(x) = x^{15}$,   $x \in [0, 1]$.

\begin{align} \label{eq:ex_1}
\begin{split}
     &{}_0^{RL} D_x^{(k)} f (x) = \frac{\Gamma(0.5)\Gamma(16)}{2\Gamma(0.5)\Gamma(15.5)} x^{14.5} +  \frac{\Gamma(1.5)\Gamma(16)}{4\Gamma(0.5)\Gamma(16.5)} x^{15.5} +  \frac{\Gamma(2.5)\Gamma(16)}{4\Gamma(0.5)\Gamma(17.5)} x^{16.5}  \\
     &f(0) = 0 , ~ f(1) = 1.
\end{split}
\end{align}

(b) Now, we solve the problem (\ref{eq:ex_2}), which has the analytical solution as $f(x) = x^{15.5}$,  $x \in [0, 1]$.

\begin{align} \label{eq:ex_2}
\begin{split}
     &{}_0^{RL} D_x^{(k)} f (x) = \frac{\Gamma(0.5)\Gamma(16.5)}{2\Gamma(0.5)\Gamma(16)} x^{15} +  \frac{\Gamma(1.5)\Gamma(16.5)}{4\Gamma(0.5)\Gamma(17)} x^{16} +  \frac{\Gamma(2.5)\Gamma(16.5)}{4\Gamma(0.5)\Gamma(18)} x^{17} , \\
     &f(0) = 0 , ~ f(1) = 1.
\end{split}
\end{align}

We seek an approximation of type (\ref{eq:jac_galerkin}), where the residual (\ref{eq:resi}) is $L^2$ orthogonal to test functions belonging to space $V_{N-1}$. The trial space $U_N$ is constructed for $\alpha = \beta = 0$ corresponding to the shifted Legendre polynomial and the test space $V_{N-1}$ is constructed using the shifted Legendre polynomial. 

Furthermore, we solve an additional equation (\ref{eq:bc}) to impose the boundary condition at $x =1$ using the Tau approach (Ch.3 \cite{canuto2006spectral}). Note that the boundary condition at $x =0$ is satisfied by the construction of trial space ($U_N \subset U$).  

\begin{equation} \label{eq:bc}
     \sum_{n=0}^N \hat{f}_n \phi_n (1) = 1
\end{equation}

\begin{figure} [!htb]
    \centering
\includegraphics[width=0.8\linewidth]{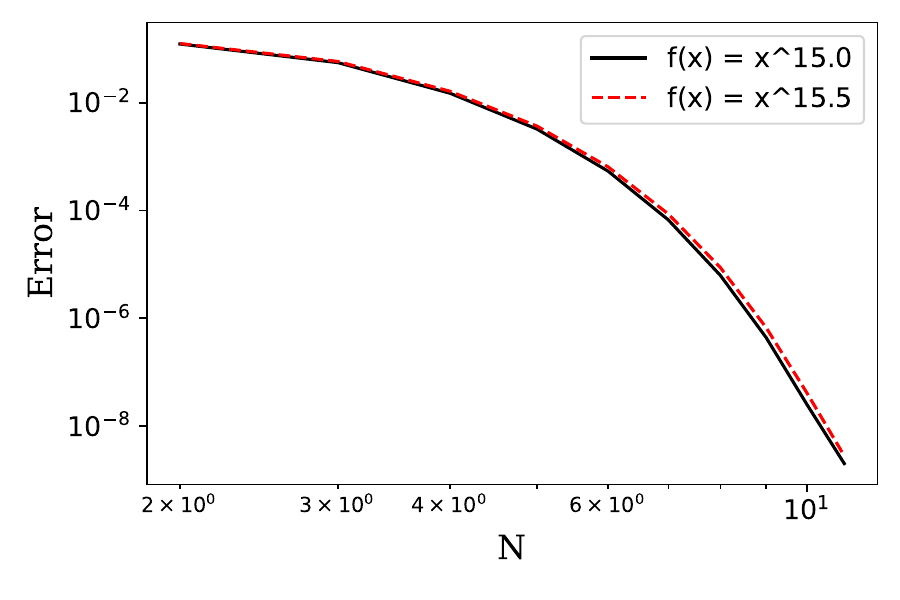}
    \caption{The rate of convergence with respect to $N$ for our Petrov-Galerkin scheme for the functions  (a) $f(x) = x^{15}$ and (b) $f(x) = x^{15.5}$. Our method converges spectrally for both polynomial and non-polynomial function.}
\label{fig:err_gen_frac_directlet}
\end{figure}

Note that, the function $f(x) = x^{15.5}$ has a non-polynomial behavior. In general, for such non-polynomial functions, approximations using polynomial leads to slow convergence (infer \cite{zayernouri2013fractional}) and vice-versa. However, we show in fig.~\ref{fig:err_gen_frac_directlet} (and table \ref{tb:err}), that our method convergences spectrally for either functions.

\subsection{Convergence analysis}

We shall now study the convergence, where we would like to bound $ \left \| f - f_N \right \|_{L^2_w [0,1]} $ and study the decay of coefficients in terms of derivative or the source term.

\begin{theorem} \label{th:decay} (Decay of coefficients) For $f \in U$, $f_N \in U_N$ and $v_N \in V_N$, then the coefficients of the Petrov-Galerkn scheme decays as,

   \begin{equation} \label{eq:decay_rate}
    |\hat{f}_n|^2  ~\leq~  \frac{C}{\lambda_n^{2m}}  \left \| g  \right \|^2_{H_w^{2m} [0,1]} ,
\end{equation}

As $m \to \infty$, the coefficients decay spectrally. 

\end{theorem}

\begin{table}[htb]
\centering
\begin{tabular}{ p{0.3\textwidth} p{0.3\textwidth} p{0.3\textwidth} }  
 \hline
$N$  & $f(x) = x^{15}$ &$f(x) = x^{15.5}$\\  
 \hline
 2 & 0.1230 & 0.1258 \\ 
4 & 0.01502 & 0.01633 \\
6 & 0.00053 & 0.000646 \\
8 & $6.3095 \times 10^{-6}$ & $8.8172 \times 10^{-6}$ \\
10 & $2.4774 \times 10^{-8}$ & $4.0516\times 10^{-8}$\\
\hline
\end{tabular} 
\caption{Mean squared error with respect to $N$ for our Petrov-Galerkin scheme for the functions  (a) $f(x) = x^{15}$ and (b) $f(x) = x^{15.5}$.} \label{tb:err}
\end{table}

\noindent
\textbf{Proof of theorem \ref{th:decay}:} Let $\phi_n = \kappa * \Tilde{P}_n^{\alpha, \beta}$, where $\Tilde{P}_n^{\alpha, \beta}$ is a nth order shifted Jacobi polynomial. The Galerkin projection of the function is given by (\ref{eq:jac_galerkin}). 

Following our construction (\ref{eq:PG_const}), we have,

\begin{equation} \label{eq:an}
   \hat{f}_n ~=~ \frac{1}{\Tilde{\gamma}_n} \left(g, \Tilde{P}_n^{\alpha, \beta} \right)_{L^2_w [0,1]}
\end{equation}

where, $ \Tilde{\gamma}_{n} = \| \Tilde{P}_n^{\alpha, \beta} \|^2_{L^2_w [0,1]}$ is the orthogonality constant for shifted Jacobi polynomials. Recall, $\Tilde{P}_n^{\alpha, \beta}$ solves the integer-order Sturm-Louville problem (\ref{eq:int_jac_stl}). 

\begin{equation} \label{eq:int_jac_stl}
    \left(\mathcal{L} + \lambda_{n} w(x) \right) \Tilde{P}_n^{\alpha, \beta} ~=~ 0 
\end{equation}

where, the weight, $w(x) = (2-2x)^{\alpha} (2x)^{\beta}$ for the shifted Jacobi polynomial ($\Tilde{P}_n^{\alpha, \beta} $), $\lambda_n = n(n+\alpha + \beta +1)$ is the corresponding $n^{th}$ eigenvalue and the differential operator, $\mathcal{L}$ is given by (\ref{eq:L_diff}).

\begin{equation} \label{eq:L_diff}
    \mathcal{L} \Tilde{P}_n^{\alpha, \beta} ~=~ \frac{d}{dx} \left ( (2-2x)^{\alpha+1} (2x)^{\beta+1} \frac{d}{dx} \Tilde{P}_n^{\alpha, \beta} \right ); ~~ \alpha, \beta > -1
\end{equation}

Using (\ref{eq:int_jac_stl}) in (\ref{eq:an}), we have, 

\begin{equation}
   \hat{f}_n ~=~ \frac{1}{\Tilde{\gamma}_n} \left(g, \Tilde{P}_n^{\alpha, \beta} \right)_{L^2_w [0,1]} ~=~ \frac{-1}{\Tilde{\gamma}_n \lambda_n} \int_0^1 g(x) \mathcal{L}  \Tilde{P}_n^{\alpha, \beta}(x) dx
\end{equation}

By performing integration by parts, we have,

\begin{equation}
   \hat{f}_n ~=~ \frac{-1}{\Tilde{\gamma}_n \lambda_n} \int_0^1 \mathcal{L} g(x)  \Tilde{P}_n^{\alpha, \beta}(x) dx
\end{equation}

We introduce the symbol $(.)_{(m)}$ defined as, 

\begin{equation}
    g_{(m)} ~=~ \frac{1}{w(x)} \mathcal{L}  g_{(m-1)} ~=~ \left(\frac{\mathcal{L}}{w(x)}\right)^m g(x),
\end{equation}

and by performing integration of part $m$-times; we have, 

\begin{equation}
   \hat{f}_n ~=~ \frac{(-1)^m}{\Tilde{\gamma}_n (\lambda_n)^m} \int_0^1 g_{(m)}  \Tilde{P}_n^{\alpha, \beta}(x) dx ~=~  \frac{(-1)^m}{\Tilde{\gamma}_n (\lambda_n)^m} \left( g_{(m)} , \Tilde{P}_n^{\alpha, \beta}\right)_{L_w^[0,1]}
\end{equation}

Now consider, $|\hat{f}_n|^2$ and apply the Cauchy-Schwartz inequality, we get (\ref{eq:decay_rate}), where, $C$ is a constant independent of $n$.

\begin{equation} 
    |\hat{f}_n|^2 ~\leq~  \frac{C}{\lambda_n^{2m}}  \left \| g_{(m)} \right \|^2_{L_w^2 [0,1]} ~\leq~  \frac{C}{\lambda_n^{2m}}  \left \| g  \right \|^2_{H_w^{2m} [0,1]} ,
\end{equation}

As $m \to \infty$, the coefficients decay spectrally.

\begin{theorem} \label{th:rate_convg} (Truncation error) For $f \in U$, $f_N \in U_N$ and $v_N \in V_N$, the Petrov-Galerkn scheme converges as,

    \begin{equation}
        \left \| f - f_N \right \|_{L^2_w [0,1]} ~\leq~ C N^{-p}  \left \|  g \right \|_{H^{p}_w [0,1]} 
    \end{equation}
\end{theorem}

\noindent
\textbf{Proof of theorem \ref{th:rate_convg}:} Consider, 

\begin{equation}
    f - \sum_{n=0}^{N} \hat{f}_n \phi_n ~=~ \sum_{n=N+1}^{\infty} \hat{f}_n \phi_n 
\end{equation}

Taking the norm and squaring it on both sides, with further use of Cauchy-Schwartz inequality to simplify of right hand side, we get

\begin{equation} \label{eq:truc_0}
    \left \| f - \sum_{n=0}^{N} \hat{f}_n \phi_n \right \|^2_{L^2_w [0, 1]} ~\leq~  \left \| \kappa \right \|^2_{L^2_w [0, 1]} \sum_{n=N+1}^{\infty} |\hat{f}_n|^2 \Tilde{\gamma}_n 
\end{equation}

where, $\Tilde{\gamma}_{n} =  \left \| \Tilde{P}_n^{\alpha, \beta} \right \|^2_{L^2_w [0,1]}$ is the orthogonality constant for shifted Jacobi polynomials. Using (\ref{eq:decay_rate}) in (\ref{eq:truc_0}), we get.

\begin{align}
    \begin{split}
         \left \| f - \sum_{n=0}^{N} \hat{f}_n \phi_n \right \|^2_{L^2_w [0, 1]} ~&\leq ~  \sum_{n=N+1}^{\infty} \frac{C  \gamma_n }{\lambda_n^{2m}} \left \| \kappa \right \|^2_{L^2_w [0, 1]}   \left \| g \right \|^2_{H_w^{2m} [0,1]} \\ 
         &\leq ~  C  N^{-4m}  \left \| \kappa  \right \|^2_{L^2_w [0, 1]}  \left \| g \right \|^2_{H_w^{2m} [0,1]} 
    \end{split}
\end{align}

\noindent
Note that, $\kappa \in L_w^2[0,1]$, hence its norm is a constant, independent of $N$. Taking square-root and choosing  $p = 2m$ and plugging, $ g = {}_0^{RL} D^{(k)}_x f $ completes the proof.

\section{Summary}

Engineering problems for real world applications are often defined over a finite domain. In this view, we first extend the results of  general fractional calculus by \citet{luchko_finite_interval} on finite interval to arbitrary orders, by introducing the Luchko class of kernels defined in (\ref{eq:luchko_class}). 

Alongside, the work of \citet{luchko2021general} for semi-infinite domains and our present work, the mathematical theory of general fractional calculus is now complete. This provides the mathematical foundations for physicists and engineers to develop mathematical models with operators of arbitrary kernels. 

Inorder to solve for general fractional differential equations, we introduced the Jacobi convolution series (\ref{eq:jconv}) (\ref{eq:jrconv}) as a first step towards development of spectral methods. It verifies they are basis functions. A notable property of this basis functions, the general fractional derivative of Jacobi convolution series is a shifted Jacobi polynomial. 

By virtue of this new class of of basis functions, we constructed a Petrov-Galerkin scheme for general fractional operators. With regards to the computational efficiency, our scheme leads to a diagonal stiffness matrix. Indeed, our approach is valid for fractional operators too, since they are a special case of general fractional operators. Our results shows that, the convergence for both polynomial and non-polynomial functions alike, which is major improvement. Following the error estimate, it is evident that, introducing such a basis function leads to methods, where the convergence rate is spectral. 

It is to be noted that, the idea of obtaining Jacobi convolution series can be extended to any arbitrary convolution type operator to develop an accurate and efficient Petrov-Galerkin scheme as long as (a) it's inverse exists and (b) it's a basis function.

\appendix

\section*{Author contributions: CRediT}
Pavan Pranjivan Mehta : Conceptualization, Methodology, Formal analysis, Investigation, Software, Validation, Project administration and Writing – original draft. Gianluigi Rozza : Funding acquisition and Resources. 

\section*{Declaration of interest}
The authors declare no conflict of interests.

\section*{Funding}
This work is funded by PNRR digitalisation, Itlay. P.I. : Gianluigi Rozza, SISSA, Italy

\bibliographystyle{elsarticle-num-names} 
\bibliography{ref}

\end{document}